\documentclass[12pt]{article}
\usepackage{mathrsfs}
\usepackage{amsmath}
\usepackage{amsfonts}
\usepackage{amssymb}
\def\endpf{\hbox{\vrule height1.5ex width.5em}}

\def\a{\alpha}
\def\b{\beta}
\newcommand{\s}{\mathcal}
\newcommand{\sm}{\setminus}
\newcommand{\CC}{{\mathbb C}}

\newcommand{\BB}{{\mathbb B}}
\newcommand{\PP}{{\mathbb P}}

\newcommand{\ra}{\rightarrow}
\newcommand{\ov}{\overline}
\def \-{\overline}
\newcommand{\p}{\partial}
\newcommand{\w}{\widetilde}
\newcommand{\wt}{\widetilde}

\newtheorem{theorem}{Theorem}[section]

\newtheorem{proposition}[theorem]{Proposition}
\newtheorem{question}[theorem]{Question}

\date{}

\begin{document}

\title{\bf  Chern-Moser-Weyl  Tensor and   Embeddings
into Hyperquadrics}


\medskip
\medskip

\author {Xiaojun Huang\footnote{ Supported in part by NSF-1363418}\
\
and Ming Xiao}

\medskip
\medskip

\vspace{3cm} \maketitle

\medskip \centerline {\it Dedicated to our friend Dick Wheeden}

\bigskip
\section{Introduction}
A central problem in Mathematics is the classification problem.
Given a set of objects and an equivalence relation, loosely
speaking, the problem asks how to find an accessible way to tell
whether two objects are in the same equivalence class. A general
approach to this problem is to find a complete set of (geometric,
analytic or algebraic) invariants. In the subject of Several Complex
Variables and Complex Geometry, a fundamental problem is to classify
complex manifolds or more generally, normal complex spaces under the
action of biholomorphic transformations. When the normal complex
spaces are open and have strongly pseudo-convex boundary, by the
Fefferman-Bochner theorem, one needs only to classify the
corresponding boundary strongly pseudoconvex CR manifolds under the
application of CR diffeomorphisms. The celebrated Chern-Moser theory
is a theory which gives two different constructions of a complete
set of invariants for such a classification problem. Among various
aspects of the Chern-Moser theory (especially the geometric aspect
of the theory), the Chern-Moser-Weyl tensor plays a key role.
However, this trace-free tensor is defined in a very complicated
manner. This makes it hard to apply in the applications. The
majority  of first several sections in this article surveys   some
work done in papers of Chern-Moser [CM], Huang-Zhang [HZh],
Huang-Zaitsev [HZa]. Here, we give a simple and more accessible
account on the Chern-Moser-Weyl tensor. We also make an immediate
application of the monotonicity property for this tensor to the
study of CR embedding problem for the positive signature case.

In the last section  of this paper, we present new materials. We
will show that the family of compact strongly pseudo-convex
algebraic hypersurfaces constructed in [HLX] cannot be locally
holomorohically embedded into a sphere of any dimension. The
argument is based on the rationality result established in [HLX] and
the Segre geometry associated with such a family.  This gives a
negative answer to a long standing folklore conjecture concerning
the embeddability of compact strongly pseudo-convex algebraic
hypersurfaces into a sphere of sufficiently high dimension. For an
extensive discussion on the history on the CR embeddability into
spheres, we refer the reader to the introduction section  of a
recent joint paper of the first author with Zaistev [HZa].

\bigskip


\section{Chern-Moser-Weyl tensor for a Levi non-degenerate hypersurface }
In this article, we assume that the CR manifolds under consideration
are already embedded as
 hypersurfaces in the complex Euclidean spaces.
 We first consider
the  case where the manifolds are even Levi non-degenerate.

We use $(z,w)\in \mathbb{C}^n\ \times\mathbb{C}$ for the coordinates
of $\mathbb{C}^{n+1}$. We always assume that $n\ge 2$, for otherwise
the Chern-Moser-Weyl tensor is identically zero. In that setting,
one has to consider the Cartan curvature functions instead, which we
will not touch  in this  article.

Let $M$ be a smooth real hypersurface.
We say that $M$ is Levi non-degenerate at $p\in M$ with signature
$\ell\le n/2$ if there is a local holomorphic change of coordinates,
that maps $p$ to the origin,  such that in the new coordinates, $M$
is defined near $0$ by an equation of the form:
\begin{equation}
r=v-|z|^2_\ell+o(|z|^2+|zu|)=0 \label{001}
\end{equation}
Here, we write $u=\Re w, v=\Im w$ and $<a,\bar b>_\ell=-\sum_{j\le
\ell} a_j \bar b_j+\sum_{j=\ell+1}^n a_j \bar b_j,
|z|_\ell^2=<z,\bar z>_\ell.$ When $\ell=0$, we regard $\sum_{j\le
\ell} a_j=0$.

Assume that $M$ is Levi non-degenerate with the same signature
$\ell$ at any point in $M$. For a point $p\in M$, a real
non-vanishing 1-form $\theta_p$ at $p\in M$ is said to be
appropriate contact form at $p$ if $\theta_p$ annihilates
$T_p^{(1,0)}+ T_p^{(0,1)}M$ and the Levi form $L_{\theta_p}$
associated with $\theta_p$ at $p\in M$ has $\ell$ negative
eigenvalues and $n-\ell$ positive eigenvalues. Here we recall the
definition of  the Levi-form $L_{\theta_p}$ at $p$ as follows: We
first extend $\theta_p$ to a smooth 1-form $\theta$ near $p$ such
that $\theta|_q$ annihilates $T_q^{(1,0)}+ T_q^{(0,1)}M$ at any
point $q\approx p$. For any $ X_{\alpha}, X_{\beta} \in T_p^{(1,0)}$, we define
\begin{equation}
 L_{\theta_p}(X_\alpha,X_\beta):=-i<d\theta|_p,
X_{\alpha}\wedge\overline X_\beta>. \label{eqn:000}
\end{equation}
One can easily verify that $L_{\theta_p}$ is a well-defined
Hermitian form in the tangent space of type $(1,0)$ of $M$ at $p$,
which is independent of the choice of the extension of the  1-form
$\theta$. In the literature, any smooth non-vanishing 1-form
$\theta$ along $M$ is called a smooth contact form,  if $\theta|_q$
annihilates $T_q^{(1,0)}M$ for any $q\in M$. If $\theta|_q$ is
appropriate at $q\in M$, we call $\theta$ an appropriate smooth
contact 1-form along $M$. Write $E_p$ for the set of appropriate
contact 1-forms at $p$ defined above, and $E$ for the disjoint union
of $E_p$. Then two elements in $E_p$ are proportional by a positive
constant for the case of $\ell<n/2$; and are proportional by a non
zero constant when $\ell=n/2.$ There is a natural smooth structure
over $E$ which makes $E$ into a $R^{+}$ fiber bundle over $M$ when
$\ell<n/2$, or a $R^{*}$-bundle over $M$ when $\ell=n/2$. When $M$
is defined  near $0$ by an equation of the form as in (\ref{001}),
then $i\partial r$ is an appropriate contact form of $M$ near $0$.
In particular, for any appropriate contact 1-form $\theta_0$ at
$0\in M$, there is a constant $c\not =0$ such that
$\theta_0=ic\partial r|_0.$ And $c>0$ when $\ell<n/2$.  Applying
further a holomorphic change of coordinates $(z,w)\ra (
\sqrt{|c|}z,cw)$ and the permutation transformation
$(z_1,\cdots,z_n,w)\ra (z_n,\cdots, z_1,w)$ if necessary, we can
simply have $\theta_0=i \partial r|_0.$ Assign the weight of
$z,\-{z}$ to be $1$ and that of $u,v,w$ to be $2$. We say
$h(z,\-{z},u)=o_{wt}(k)$ if $\frac{h(tz,\-{tz},t^2u)}{t^k}\ra 0$
uniformly on compact sets in $(z,u)$ near the origin. We write
$h^{(k)}(z,w)$ for a weighted homogeneous holomorphic polynomial of
weighted degree $k$ and $h^{(k)}(z,\-{z},u)$ for a weighted
homogeneous polynomial of weighted degree $k$. We first have the
following  special but crucial case of the Chern-Moser normalization
theorem:

\bigskip
{\proposition {\it Let $M\subset {\mathbb C}^n\times {\mathbb C}$ be a
smooth Levi non-degenerate hypersurface. Let $\theta_p\in E_p$ be an
appropriate real 1-form at $p\in M$. Then there is a biholomorphic
map $F$ from a neighborhood of $p$ to a neighborhood of $0$   such
that $F(p)=0$ and  $F(M)$  near $0$ is defined by an equation of the
following normal form (up to fourth order):
\begin{equation}
r=v-|z|_\ell^2+\frac{1}{4}s
(z,\bar{z})+R(z,\-{z},u)=v-|z|_\ell^2+\frac{1}{4}\sum s^0_{\alpha\bar
{\beta}\gamma \bar{\delta}}z_{\alpha} {\bar z_\beta}z_{\gamma}{\bar
z_\delta}+R(z,\-{z},u)=0. \label{eqn:002}
\end{equation}

Here  $s(z,\-{z})=\sum s^0_{\alpha\bar \beta
\gamma\bar{\delta}}z_{\alpha} {\bar z_\beta}z_{\gamma}{\bar
z_\delta}$,$\ s^0_{\alpha\bar {\beta}\gamma \bar{\delta}}=
s^0_{\gamma \bar {\beta}\alpha\bar{\delta}}= s^0_{\gamma
\bar{\delta}\alpha\bar {\beta}},\ \overline {s^0_{\alpha\bar
{\beta}\gamma
\bar{\delta}}}=s^0_{\beta\bar{\alpha}\delta\bar{\gamma}}$
 and
\begin{equation}
 \sum_{\alpha, \beta=1}^n s^0_{\alpha\bar {\beta}\gamma
\bar{\delta}}g_0^{\bar \beta \alpha}=0 \label{eqn:002-1}
\end{equation}
where $g_0^{\bar \beta \alpha}=0$ for $\beta\neq\alpha$,
$g_0^{\bar\beta\beta}=1$ for $\beta>\ell, g_0^{\bar\beta\beta}=-1$
for $\beta\leq \ell$.  Also $R(z,\-{z},u)=o_{wt}(|(z,u)|^4)\cap
o(|(z,u)|^4)$. Moreover, we have  $i\p r|_0=(F^{-1})^*\theta_p.$}}

\bigskip
{\it Proof of Proposition 2.1}: By what we discussed above, we can assume
that $p=0$ and $M$ near $p=0$ is defined by an equation of the form
as in (\ref{001}).  We first show that we can get rid of all
weighted third order degree terms. For this purpose, we choose a
transformation of the form $f=z+f^{(2)}(z,w)$ and
$g=w+g^{(3)}(z,w)$. Suppose that $F=(f_1,\cdots,f_n,g)=(f,g)$ maps
$(M,p=0)$ to a hypersurface near $0$ defined by an equation of the
form as in (\ref{001}) but without weighted degree $3$ terms in the
right hand side. Substituting $F$ into the new equation and
comparing terms of weighted degree three, we get
$$\Im \left( g^{(3)}-2i<\-{z},f^{(2)}>_\ell\right)|_{w=u+i|z|_\ell}=G^{(3)}(z,\-{z},u)$$
where $G^{(3)}$ is a certain given real-valued polynomial of
weighted degree $3$ in $(z,\-{z},u)$. Write
$G^{(3)}(z,\-{z},u)=\Im\{a^{(1)}(z)w+\sum_{j=1}^{n}b^{(2)}_{j}(z)\-{z_j}\}.$
Choosing  $g^{(3)}=a^{(1)}(z)w$ and
$f_j^{(2)}=\frac{i}{2}b_j^{(2)}(z),$  it then  does our job.

Next, we choose a holomorphic transformation of the form
$f=z+f^{(3)}(z,w)$ and
$g=w+g^{(4)}(z,w)$ to simplify the weighted degree $4$ terms in the defining equation of $(M,p=0)$.
Suppose that $M$ is originally defined by
$$r=v-|z|^2_\ell+A^{(4)}(z,\-{z},u)+o_{wt}(4)=0 \label{001-n}$$
and is transformed to an equation of the form:
$$r=v-|z|^2_\ell+N^{(4)}(z,\-{z},u)+o_{wt}(4)=0 \label{001-n}.$$
substituting the map $F$ and collecting terms of weighted degree $4$, we get the equation:
$$\Im \left( g^{(4)}-2i<\-{z},f^{(3)}>_\ell\right)|_{w=u+i|z|_\ell}=N^{(4)}(z,\-{z},u)-A^{(4)}(z,\-{z},u).$$
Now, we like to make $N^{(4)}$ as simple as possible by choosing $F$.
Write
$$-A^{(4)}=\Im\{ b^{(4)}(z)+b^{(2)}(z)u+b^{(0)}u^2+\sum_{j=1}^{n}c_j^{(3)}(z)\-{z_j}+\sum_{|\a|=|\b|=2}\wt{c_{\a\-{\b}}}z^\a\-{z^\b}\}.$$
Let $$X^{(4)}(z,w)=b^{(4)}(z)+b^{(2)}(z)w+b^{(0)}w^2,\
-2i\delta_{j\ell}Y_j^{(3)}(z,w)=c_j^{(3)}(z)-
ib^{(2)}(z)z_j-2ib^{(0)}z_jw,$$
$$ Y^{(3)}=(Y_1^{(3)},\cdots, Y_n^{(3)}),$$
where $\delta_{j\ell}$ is $1$ for $j>\ell$ and is $-1$ otherwise.
Then $\Im \left(
Y^{(4)}-2i<\-{z},X^{(3)}>_\ell\right)+A^{(4)}(z,\-{z},u)=-\Im({b^{(0)}})|z|_\ell^4+\sum_{|\a|=|\b|=2}d_{\a\-{\b}}z^\a\-{z^\b}.$
By the Fischer decomposition theorem ([SW]), write in the  unique
way
$$-\Im({b^{(0)}})|z|_\ell^4+\sum_{|\a|=|\b|=2}d_{\a\-{\b}}z^\a\-{z^\b}=h^{(2)}(z,\-{z})|z|_{\ell}+h^{(4)}(z,\-{z}).$$
Here $h^{(2)}(z,\-{z})$ and $h^{(4)}(z,\-{z})$ are real-valued,
bi-homogeneous in $(z,\overline{z})$ and $\Delta_\ell
h^{(4)}(z,\-{z})=0$. Here, we write $\triangle_\ell=-\sum_{j\le
\ell} \frac{\p^2}{\p z_j\p\bar z_j} +\sum_{j=\ell+1}^n
\frac{\p^2}{\p z_j\p\bar z_j}$. Notice that $h^{(2)}$ has no
harmonic terms, we can find $Z^{(1)}(z)$ such that $ \Re(<\-{z},
Z^{(1)}(z)>_\ell)=0$ and $\Im(2<\-{z},Z^{(1)}>)=h^{(2)}(z,\-{z}).$
Finally, if we define $f=z+X^{(4)}(z,w)+Z^{(1)}(z)w$ and
$g^{(4)}=w+Y^{(4)}$, then $(f,g)$ maps $(M,0)$ to a hypersurface
with $R(z,\-{z},u)=o_{wt}(4)\cap O(|(z,u)|^3)$. Now suppose that the
terms with non-weighted degree of $3 $ or $4$ in $R$ are uniquely
written as
$ub^{(3)}(z,\ov{z})+u^2\Im{(b^{(1)}(z))+b^{(0)}u^3+c^{(0)}u^4}$ with
$b^{(3)}(z,\ov{z})=\Im{(c^{(3)}(z)+\sum_{|\a|=2,|\b|=1}d_{\a\ov{\b}}z^\a\ov
z^\b)}.$  Then we need to make further change of variables as
follows to make $R=o_{wt}(4)\cap o(|(z,u)|^4)$ without changing
$N^{(4)}(z,\-{z})$:
$$w'=w+wc^{(3)}(z)+w^2b^{(1)}(z)+ib^{(0)}w^3+ic^{(0)}w^4,$$
$$z'_j=z_j+\delta_{j,\ell}wb^{(1)}(z)z_j+\frac{i}{2}\sum_{|\a|=2}wd_{\a,\overline{j}}z^\a+\delta_{j,\ell}\frac{3i}{2}w^2z_jb^{(0)}.$$

 Now, the
trace-free condition in (\ref{eqn:002-1}) is equivalent to the
following condition :
\begin{equation*}
\triangle_\ell s(z,\bar z)\equiv 0.
\end{equation*}
Indeed, this follows from the following  fact: Let $
\Delta_{H}=\sum_{l,k=1}^{n}h^{l\-{k}}\partial_l\-{\partial}_k$ with
$\-{h^{l\-{k}}}=h^{k\-{l}}$ for any $l,k$. Then
\begin{equation}
\Delta_{H}s^0 (z,\-{z})=4\sum_{\gamma,\delta=1}^{n}\sum_{\a,
\b=1}^{n}h^{\a\-{\b}}s^0_{\a\-{\b}\gamma\-{\delta}}z_\gamma\-{z_\delta}.
\end{equation}
 This proves the proposition. $\endpf$
\bigskip


We assume the notation and conclusion in Proposition 2.1.  The
Chern-Moser-Weyl tensor at $p$ associated with the appropriate
1-form $\theta_p$ is defined as
the 4th order  tensor  $S_{\theta_p}$
acting over
$T_p^{(1,0)}M\otimes T_p^{(0,1)}M\otimes T_p^{(1,0)}M\otimes
T_p^{(0,1)}M$. More precisely,   for each $X_p,Y_p.Z_p,W_p\in
T_p^{(1,0)}M$, we have the following definition:

Let $F$ be the biholomorphic map sending $M$ near $p$ to the normal
form as in Proposition 2.1 with $F(p)=0$, and
 write
$F_*(X_p)=\sum_{j=1}^{n}a^j\frac{\partial}{\partial z_j}|_0:=X_p^0$,
$F_*(Y_p)=\sum_{j=1}^{n}b^j\frac{\partial}{\partial z_j}|_0:=Y_p^0,$
$F_*(Z_p)=\sum_{j=1}^{n}c^j\frac{\partial}{\partial z_j}|_0:=Z_p^0,$
and $F_*(W_p)=\sum_{j=1}^{n}d^j\frac{\partial}{\partial
z_j}|_0:=W_p^0.$ Then

\begin{equation}
S_{\theta_p}(X_p,\-{Y_p},Z_p,\-{W_p}):=
\sum_{\alpha,\beta,\gamma,\delta=1}^{n}s^0_{\alpha\-{\beta}\gamma\-{\delta}}a^\alpha\-{b^\beta}c^{\gamma}
\-{d^\delta},\ \ \hbox{which is denoted by}~
S_{\theta_0}(X^0_p,\-{Y^0_p},Z^0_p,\-{W^0_p}). \label{eqn:002-02}
\end{equation}

Since the normalization map $F$ is not unique, we have to verify
that the tensor $S_{\theta_p}$ is well-defined. Namely, we need to
show that it is independent of the choice of the normal coordinates.
We  do this in the next section. For the rest of this section, we
assume this fact and derive some basic properties for the tensor.

 For a  basis
$\{X_{\alpha}\}_{\alpha=1}^n$ of $ T_p^{(1,0)}M$ with $p\in M$,
write $({ S_{{\theta}_p}})_{\alpha\bar {\beta}\gamma
\bar{\delta}}= S_{\theta_p}(X_{\alpha}, \overline X_{\beta},
X_{\gamma}, \overline X_{\delta})$. From the definition, we then
have the following symmetric properties:
\begin{eqnarray*}
&(S_{{\theta}_p})_{\alpha\bar {\beta}\gamma
\bar{\delta}}=(S_{\theta_p})_{\gamma\bar {\beta}\alpha
\bar{\delta}}=(S_{\theta_p})_{\gamma \bar{\delta}\alpha\bar
{\beta}}\\
&\overline{(S_{\theta_p})_{\alpha \bar{\beta}\gamma \bar{\delta}}}
=(S_{\theta_p})_{\beta\bar{\alpha} \delta\bar{\gamma}},
\end{eqnarray*}
and the following trace-free condition:
\begin{equation}
\sum_{\beta,\alpha=1}^n g^{\bar \beta
\alpha}(S_{\theta_p})_{\alpha\bar {\beta}\gamma
\bar{\delta}}=0.\label{eqn:000-05}
\end{equation}
Here
\begin{equation}
 g_{\alpha\bar\beta}=L_{\theta|_p}(X_\alpha,X_\beta):=-i<(d\theta)|_p,
X_{\alpha}\wedge\overline X_\beta> \label{eqn:000-01}
\end{equation}
is the Levi form of $M$ associated with $\theta_p$ and $\theta$ is a
smooth extension of $\theta_p$ as a proper contact form of $M$ near
$p$. Also, $(g^{\bar \beta \alpha})$ is the inverse matrix of
$(g_{\alpha\bar \beta}).$ In the following, we write
$\wt{\theta}=(F^{-1})^*(\theta).$

To see the trace-free property in (\ref{eqn:000-05}), we write that
$F_*(X_\alpha)=\sum_{k=1}^{n}a_\alpha^k\frac{\partial}{\partial
z_k}|_0.$ Then
$g_{\alpha\bar\beta}=L_{\theta_p}(X_\alpha,X_\beta)=-i<(d\theta)|_p,
X_{\alpha}\wedge\overline X_\beta>=-i<(dF^*(\wt{\theta})|_p,
X_{\alpha}\wedge\overline X_\beta>=-i<(i\-{\partial}{\partial}r|_0,
F_*(X_{\alpha})\wedge\overline
{F_*(X_\beta)}>=(g_0)_{k\-{l}}a_\alpha^k\-{a_\beta^l}.$ Here
$(g_0)_{k\-{l}}$ is defined as before. Write $G=(g_{\a\b}),
G^0=(g_0)_{\a\b}, A=(a_k^l), B=A^{-1}:=(b_k^l)$. Then we have the
matrix relation: $G=AG^0\-{A}^t$. Thus
$G^{-1}=(\-{A^t})^{-1}(G^0)^{-1}A^{-1},$ from which we have
$g^{\gamma \-{\b}}=\-{b^{\b}_{l}}(g_0)^{j\-{l}}b_j^\gamma.$ Thus,

$$g^{\a\-{\b}}S_{\a\-{\b}\gamma\-{\delta}}=\-{b^{\b}_{l}}(g_0)^{j\-{l}}b_j^\a
s^0_{\wt{k}\-{\wt{j}}\wt{l}\-{\wt{m}}}a_\alpha^{\wt{k}}\-{a_\beta^{\wt{j}}}a_{\gamma}^{\wt{l}}
\-{a_\delta^{\wt{m}}}=(g_0)^{j\-{l}}s^0_{j\-{l}\wt{l}\-{\wt{m}}}a_{\gamma}^{\wt{l}}
\-{a_\delta^{\wt{m}}}=0.$$


We should mention the above argument can also be easily adapted to
show the biholomorphic invariance of the appropriateness. Namely, if
$F$ is a CR diffeomorphism between two Levi non-degenerate
hypersurfaces  $M$ and $\wt{M}$ of signature $\ell$. For
$\wt{\theta_q}$ is an appropriate contact 1-form at $q\in \wt{M}$,
then $F^*(\wt{\theta_q})$ is also an appropriate contact 1-form at
$F^{-1}(q)\in M$.

For a smooth vector field $X, Y, Z, W$ of type $(1,0)$ and an
appropriate smooth contact form along $M$, ${\cal
S}_{\theta}(X,\-{Y},Z,\-{W})$ is also a smooth function along $M$.
One easy way to see this is to use the Webster-Chern-Moser-Weyl
formula obtained in [We1] through the curvature tensor of the
Webster pseudo-Hermitian metric, whose constructions are done  by
only applying the algebraic and differentiation operations on the
defining function of $M$. Another more direct way is to trace the
dependence of the tensor on the base points under the above
normalization procedure.

\medskip
Assume that $\ell>0$ and define
\begin{equation*}
\mathcal{C}_\ell=\{z\in\CC^n: |z|_\ell=0\}.
\end{equation*}
Then $\mathcal C_\ell$ is a real algebraic variety of real
codimension $1$ in $\CC^n$ with the only singularity at $0$.  For
each $p\in M$, write $\s C_\ell T_p^{(1,0)}M=\{v_p\in T_p^{(1,0)}M:\
<(d\theta)|_p, v_p\wedge{\bar v_p}>=0\}.$  Apparently, $\s
C_lT_p^{(1,0)}M$ is independent of the choice of $\theta_p$. Let $F$
be a CR diffeomorphism from $M$ to $M'$. We also have $F_{*}(\s
C_\ell T_p^{(1,0)}M)=C_\ell T_{F(p)}^{(1,0)}M'$.  Write $\s C_\ell
T^{(1,0)}M=\coprod_{p\in M} \s C_\ell T_p^{(1,0)}M$ with the natural
projection $\pi$ to $M$. We say that $X$ is a smooth section of $\s
C_\ell T^{(1,0)}M$ if $X$ is a smooth vector field of type $(1,0)$
along $M$ such that $X|_p\in \s C_\ell T_p^{(1,0)}M$ for each $p\in
M$.  $\s C_\ell T^{(1,0)}M$ is a kind of smooth bundle with each
fiber isomorphic to $\s C_\ell$.

$\s C_\ell$ is obviously a uniqueness set for holomorphic functions.
The following lemma shows that it is also a uniqueness set for the
Chern-Moser-Weyl curvature tensor. (For the proof, see Lemma 2.1 of
[HZh].)

\begin{proposition} (Huang-Zhang [HZh])
(I). Suppose that $H(z,\bar z)$ is a real real-analytic function in
$(z, \bar z)$ near $0$. Assume that $\triangle_\ell H(z, \bar
z)\equiv 0$ and $H(z, \bar z)|_{\s C_\ell}=0$. Then $H(z, \bar
z)\equiv 0$ near $0$. (II). Assume the above notation and $\ell>0$.
If $ S_{\theta_p}(X,\-X, X, \- X)=0$ for any $X\in \s C_\ell
T^{(1,0)}_pM$, then $ S_{\theta|_p}\equiv 0.$
\end{proposition}

\bigskip
\section {Transformation law for the Chern-Moser-Weyl tensor}
We next show that the Chern-Moser-Weyl tensor defined in the
previous section is well-defined by proving a transformation law. We
follow the approach and expositions developed  in Huang-Zhang [HZh].

Let $\widetilde M \subset \CC^{N+1}=\{(z, w)\in \CC^n\times \CC\}$
be also a Levi non-degenerate real hypersurface near $0$ of
signature $\ell\ge 0$ defined by an equation of the form:

\begin{equation}
\w{r}=\Im \w{w}-|\w{z}|^2_\ell+o(|\w{z}|^2+|\w{z}\w{u}|)=0.
\label{002}
\end{equation}
Let $F:=(f_1, \ldots, f_n,\phi, g): M\rightarrow \widetilde M$ be a
smooth  CR diffeomorphism.
 Then, as in [Hu1] and \cite{BH}, we can write
\begin{eqnarray}\label{eqn:022}
\begin{aligned}
\tilde z&=\tilde f(z, w)=(f_1(z, w), \ldots, f_n(z, w))=\lambda
z U+\vec{a}w+ O(|(z, w)|^2)\\
\tilde w&=g(z, w)=\sigma\lambda^2w+O(|(z, w)|^2).
 \end{aligned}
\end{eqnarray}

 Here $U\in SU(n, \ell)$. (Namely $<X U, Y\overline{ U}>_\ell=<X,
Y>_\ell$ for any $X, Y\in \CC^n$). Moreover, $\ \vec{a}\in \CC^n,\
\lambda
>0$ and $\sigma=\pm1$ with $\sigma= 1$ for $\
\ell<\frac{n}{2}$. When $\sigma=-1$, by considering
$F\circ\tau_{n/2}$ instead of $F$, where
 $\tau_{\frac{n}{2}}(z_1,\ldots,
z_{\frac{n}{2}},z_{\frac{n}{2}+1},\ldots,
 z_n, w)=(z_{\frac{n}{2}+1},\ldots,
 z_n,z_1,\ldots,
z_{\frac{n}{2}},-w),$ we can make $\sigma=1$. Hence, we will assume
in what follows that $\sigma=1$.

Write $r_0=\frac{1}{2}\Re \{g^{''}_{ww}(0)\},\ q(\tilde z, \tilde
w)=1+2i<\tilde z,  \lambda^{-2}\overline{\vec
a}>_\ell+\lambda^{-4}(r_0-i|\vec a|_\ell^2)\tilde w$,
\begin{equation}\label{eqn:033}
T(\tilde z, \tilde w)=\frac{(\lambda^{-1}(\tilde z- \lambda^{-2}\vec
a \tilde w) U^{-1}, \lambda^{-2}\tilde w)}{q(\tilde z, \tilde w)}.
\end{equation}
 Then
\begin{equation}
 F^{\sharp}(z,w)=(\tilde f^{\sharp}, g^{\sharp})(z,w): =T\circ F(z, w)=(z, w)+O(|(z,
 w)|^2)
\label{eqn:044}
\end{equation}
 with $\Re \{g^{\sharp''}_{ww}(0)\}=0$.

 Assume that $\widetilde M$  is also defined in the Chern-Moser normal
 form up to the 4th order:
 \begin{equation}
 \tilde r=\Im \tilde w-|\tilde z|_\ell^2+\frac{1}{4}\tilde {s}(\tilde
 z, \bar{\tilde z})+o_{wt}(|(\tilde z,\wt{u})|^4)=0.
\label{eqn:003}
 \end{equation}
Then $M^\sharp= T(\widetilde M)$ is defined by
\begin{equation}
r^{\sharp}=\Im w^{\sharp}-|z^{\sharp}|_\ell^2+\frac{1}{4}
s^{\sharp}(z^{\sharp},
\bar{z^{\sharp}})+o_{wt}(|(z^{\sharp},w^{\sharp})|^4)=0 \label{eqn:004}
\end{equation}
with $s^{\sharp}(z^{\sharp}, \bar{z^{\sharp}})= \lambda^{-2}\tilde
{s}(\lambda z^{\sharp}U, \lambda \overline{z^{\sharp} U})$.

One can verify that
\begin{equation}
(-\sum_{j=1}^\ell \frac{\p^2}{\p z_j^{\sharp}\p\bar z_j^{\sharp}}
+\sum_{j=\ell+1}^N \frac{\p^2}{\p z_j^{\sharp}\p\bar z_j^{\sharp}})
s^{\sharp}(z^{\sharp}, \overline {z^{\sharp}})=0.
 \label{eqn:005}
\end{equation}
Therefore (\ref{eqn:004}) is also in the Chern-Moser normal form up
to the 4th order.
 Write $F^{\sharp}(z, w)=\sum_{k=1}^{\infty}F^{\sharp
(k)}(z, w)$. Since $F^{\sharp}$ maps $M$ into
$M^{\sharp}=T(\widetilde M)$, we get the following
\begin{eqnarray}\label{22}
\begin{aligned}
 \Im\{\sum_{k\geq2}g^{\sharp(k+1)}(z,
w)-2i\sum_{k\geq2}<f^{\sharp(k)}(z, w), \bar z>_\ell\}&\\
=\sum_{k_1,\ k_2\geq2}<f^{\sharp(k_1)}(z, w), \overline{
f^{\sharp(k_2)}(z, w)}>_\ell&+\frac{1}{4}( s(z, \bar z)-
s^{\sharp}(z,\overline{z}))+o_{wt}(4)
\end{aligned}
\end{eqnarray} over $\Im w =|z|_\ell^2$.
Here, we write $F^{\sharp}(z, w)=(f^{\sharp}(z,w), g^{\sharp}(z, w))$.\\
Collecting terms of weighted degree 3 in (\ref{22}), we get
\begin{equation*}
\Im\{g^{\sharp(3)}(z, w)-2i<f^{\sharp(2)}(z, w), \bar z>_\ell\}=0\ \
\text{on} \ \ \Im w =|z|_\ell^2.
\end{equation*}
By \cite{Hu}, we get $g^{\sharp(3)}\equiv 0, f^{\sharp(2)}\equiv
0$.\\
Collecting terms of weighted degree 4 in (\ref{22}), we get
\begin{equation*}
\Im \{g^{\sharp (4)}(z,w)-2i<f^{\sharp(3)}(z, w), \bar z>_\ell\}=
\frac{1}{4}(  s(z, \bar z)- s^{\sharp}(z, \overline{z})).
\end{equation*}
Similar to the argument in \cite{Hu} and making use of the fact that
$\Re\{\frac{\p ^2 g^{\sharp(4)}}{\p w^2}(0)\}=0$, we get the
following:
\begin{eqnarray}\label{33}
\begin{aligned}
&g^{\sharp(4)}\equiv 0,\ f^{\sharp(3)}(z,w)=\frac{i}{2}a^{(1)}(z)w,
\\
<a^{(1)}(z), \bar z>_\ell|z|_\ell^2=&\frac{1}{4}( s(z, \bar z)
-s^{\sharp}(z, \overline {z}))=\frac{1}{4}( s(z, \bar z)
-\lambda^{-2}\widetilde {s}(\lambda z U,
 \overline{\lambda z U})).
\end{aligned}
\end{eqnarray}
Since the right hand side of the above equation is annihilated by
$\Delta_\ell$ and the left hand side of the above equation is
divisible by $|z|^2_\ell$. We conclude that $f^{\sharp(3)}(z,w)=0$
and

\begin{equation} \label{eqn:33-001}
 s(z, \bar z)=\lambda^{-2}\widetilde {s}(\lambda z U,
 \overline{\lambda z U}).
\end{equation}

Write $\theta_0=i\partial r|_0$ and  $\wt{\theta_0}=i\partial
\wt{r}|_0$. Then $F^*(\wt{\theta_0})=\lambda^2\theta_0.$ For any
$X=\sum_{j=1}^{n} z_j\frac{\partial}{\partial z_j}|_0,$
$F_*(X)=\lambda (z_1\frac{\partial}{\partial \wt{z_1}}|_0,\cdots,
z_n\frac{\partial}{\partial \wt{z_n}}|_0) U.$ Under this notation,
(\ref{eqn:33-01}) can be written as
$$S^0_{F^*(\wt{\theta_0})}(X,\-{X},X,\-{X})=S^0_{\wt{\theta_0}}(F_*(X),\-{F_*(X)},F_*(X),\-{F_*(X)}).$$
This immediately gives the following transformation law and thus the
following theorem, too.
\begin{equation} \label{eqn:33-01}
S^0_{F^*(\wt{\theta_0})}(X,\-{Y},Z,\-{W})=S^0_{\wt{\theta_0}}(F_*(X),\-{F_*(Y)},F_*(Z),\-{F_*(W)}),
\hbox{ for } X, Y, Z, W\in T_0^{(1,0)}M.
\end{equation}

 \bigskip
 {\theorem {\it  (1). The Chern-Moser-Weyl tensor defined in the
 previous section is independent of the choice of the  normal coordinates and thus
 is a well-defined fourth order tensor.  (2). Let $F$ be a CR diffeomorphism
 between two
  Levi non-degenerate hypersurfaces $M, M'\subset {\CC}^{n+1}$. Suppose $F(p)=q$. Then, for any appropriate contact 1-form $\wt{\theta_q}$ of $\wt{M}$ at $q$
  and a vector $v\in T_p^{(1,0)}M,$ we have the following transformation
  formula for the corresponding Chern-Moser-Weyl tensor:
  \begin{equation}
\tilde { S}_{\tilde \theta_p}(F_*(v_1), \overline {F_*(v_2)},
F_*(v_3), \overline {F_*(v_4)})=  S_{F^*(\tilde{\theta}_q)}(v_1, \-{
v_2}, v_3, \-{ v_4}). \label{eqn:303}
\end{equation}
}}

\bigskip
{\it Proof}: Let $\theta_p$ be an appropriate contact form of $M$ at
$p$, and let $F_1, F_2$ be two normalization (up to fourth order) of
$M$ at $p$. Suppose that $F_1(M)$ and $F_2(M)$ are defined near $0$
by  equations $r_1=0$ and $r_2=0$ as in (\ref{001}),
respectively. Write $\Phi=F_2\circ F_1^{-1}$ and
$\theta^1_0=i\partial r_1$, $\theta^2_0=i\partial r_2$. We also
assume that $F_1^*(\theta^1_0)=\theta_p$ and
$F_2^*(\theta^2_0)=\theta_p$. Then
 for any $X_p, Y_p,Z_p,W_p\in T_p^{(1,0)}M$,  we have
$$S^1_{\theta_p}(X_p,
\-{Y_p},Z_p,\-{W_p})=S^1_{\theta^1_0}((F_1)_{*}(X_p),
\-{(F_1)_{*}(Y_p)}, (F_1)_{*}(Z_p), \-{(F_1)_{*}(W_p)})$$ if we
define the tensor at $p$ by applying $F_2$.  We also have
$$S^2_{\theta_p}(X_p,
\-{Y_p},Z_p,\-{W_p})=S^2_{\theta^2_0}((F_2)_{*}(X_p),
\-{(F_2)_{*}(Y_p)}, (F_2)_{*}(Z_p), \-{(F_2)_{*}(W_p)}),$$ if we
define the tensor at $p$ by applying $F_2$. Since
$\theta^2_0=\Phi^*(\theta_0^1)$, and
$\Phi_*((F_1)_{*}(X_p))=(F_2)_*(X_p)$, by the transformation law
obtained in (\ref{eqn:33-01}), we see the proof in Part I of the
theorem. The proof in Part II of the theorem also follows easily
from the formula in (\ref{eqn:33-01}).

\section {A monotonicity  theorem for the Chern-Moser-Weyl tensor}
We now let $M_\ell\subset {{\mathbb C}^{n+1}} $ be a Levi
non-degenerate hypersurface with signature $\ell>0$ defined in the
normal form as in (\ref{eqn:002}).  Let $F=(f_1,\cdots,f_N,g)$ be a
CR-transversal CR embedding from $M_\ell$ into ${\mathbb
H}^{N+1}_\ell$ with $N\ge n$. Then again as in Section 3,  a simple linear algebra argument
([HZh]) shows that after a holomorphic change of variables, we can
make $F$ into the following preliminary normal form:
\begin{eqnarray}\label{eqn:022}
\begin{aligned}
\tilde z&=\tilde f(z, w)=(f_1(z, w), \ldots, f_N(z, w))=\lambda
z U+\vec{a}w+ O(|(z, w)|^2)\\
\tilde w&=g(z, w)=\sigma\lambda^2w+O(|(z, w)|^2).
 \end{aligned}
\end{eqnarray}
Here $U$ can be extended to an $N\times N$ matrix $\widetilde U\in
SU(N, \ell)$. Moreover, $\
\vec{a}\in \CC^N,\  \lambda
>0$ and $\sigma=\pm1$ with $\sigma= 1$ for $\
\ell<\frac{n}{2}$. When $\sigma=-1$, qs discussed before, by considering
$F\circ\tau_{n/2}$ instead of $F$, where
 $\tau_{\frac{n}{2}}(z_1,\ldots,
z_{\frac{n}{2}},z_{\frac{n}{2}+1},\ldots,
 z_n, w)=(z_{\frac{n}{2}+1},\ldots,
 z_n,z_1,\ldots,
z_{\frac{n}{2}},-w),$ we can make $\sigma=1$. Hence, we will assume
 that $\sigma=1$.

Write $r_0=\frac{1}{2}\Re \{g^{''}_{ww}(0)\},\ q(\tilde z, \tilde
w)=1+2i<\tilde z,  \lambda^{-2}\overline{\vec
a}>_\ell+\lambda^{-4}(r_0-i|\vec a|_\ell^2)\tilde w$,
\begin{equation}\label{eqn:033}
T(\tilde z, \tilde w)=\frac{(\lambda^{-1}(\tilde z- \lambda^{-2}\vec
a \tilde w)\widetilde U^{-1}, \lambda^{-2}\tilde w)}{q(\tilde z,
\tilde w)}.
\end{equation}
 Then
\begin{equation}
 F^{\sharp}(z,w)=(\tilde f^{\sharp}, g^{\sharp})(z,w): =T\circ F(z, w)=(z, 0, w)+O(|(z,
 w)|^2)
\label{eqn:044}
\end{equation}
 with $\Re \{g^{\sharp''}_{ww}(0)\}=0$.
Now,  $T({\mathbb H}^{N+1}_\ell)={\mathbb H}^{N+1}_\ell$. With the
same argument as in the previous section, we also arrive at the
following:

\begin{eqnarray}\label{33}
\begin{aligned}
g^{\sharp(3)}=g^{\sharp(4)}\equiv 0,\
f^{\sharp(3)}(z,w)=\frac{i}{2}a^{(1)}(z)w, &
\\
<a^{(1)}(z), \bar
z>_\ell|z|_\ell^2=|\phi^{\sharp(2)}(z)|^2+&\frac{1}{4} s(z, \bar z).
\end{aligned}
\end{eqnarray}
In the above equation, if we let $z$ be such that $|z|_\ell=0$, we
see that $s(z, \overline z)\le 0$. Now, if $F$ is not CR
transversal but not totally non-degenerate in the sense that $F$ does not map an open subset of ${\mathbb C}^n$ into ${\mathbb H}^{N}_\ell$ (see [HZh]), then one can apply this result on a dense open subset
of $M$ [BER]  where $F$ is CR transversal and then take a limit as
did in [HZh]. Then we have the following special case of the
monotonicity theorem for the Chern-Moser-Weyl tensor obtained in
Huang-Zhang [HZh]:

\begin{theorem} ([HZh]) \label {thm111} Let $M_\ell\subset \CC^{n+1}$ be a Levi
non-degenerate real hypersurface of signature $\ell$. Suppose that
$F$ is a holomorphic mapping defined in a (connected) open
neighborhood $U$ of $M$ in ${\mathbf C}^{n+1}$ that sends  $M_\ell$
into ${\mathbf H}_\ell ^{N+1}\subset \CC^{N+1}$. Assume that
$F(U)\not \subset {\mathbf H}_\ell^{N+1}$. Then when
$\ell<\frac{n}{2}$, the Chern-Moser-Weyl curvature tensor with
respect to any appropriate contact form $\theta$ is pseudo
semi-negative in the sense that for any $p\in M$, the following
holds:
\begin{equation}
 \s{S}_{ \theta|_{p}}(v_p,
\overline{v_p}, v_p, \overline {v_p})\le 0,\ \ \hbox{for}\ v_p\in
{\s C}_\ell T_p^{(1,0)}M.
\end{equation}
 When
$\ell=\frac{n}{2}$, along a certain contact form $\theta$, $\s
S_{\theta}$ is pseudo negative.
\end{theorem}

\section {Counter-examples to the embeddability problem for
compact algebraic Levi non-degenerate  hypersurfaces with positive
signature into hyperquadrics}

In this section, we apply Theorem \ref{thm111} to construct a
compact Levi-nondegenerate hypersurface in a projective space, for
which any piece of it can not be holomorphically embedded into a
hyperquadric of any dimension with the same signature. This section
is based on the work in  the last section of Huang-Zaitsev [HZa].

Let $n,\ell$ be two integers with $1<\ell\le n/2$. For any
$\epsilon$, define
$${M_\epsilon}:=
\left\{[z_0,\cdots,z_{n+1}]\in {\PP}^{n+1}:
 |z|^2 \left(-\sum_{j=0}^{\ell} |z_j|^2 + \sum_{j=\ell+1}^{n+1}|z_j|^2\right)
 +\epsilon\left(|z_1|^4-|z_{n+1}|^4\right)
=0 \right\}.$$ Here $|z|^2=\sum_{j=0}^{n+1}|z_j|^2$ as usual. For
$\epsilon =0$, ${M_\epsilon}$ reduces to the generalized sphere with
signature $\ell$, which is the boundary of the generalized ball
$${\BB}^{n+1}_\ell:=
\left\{\{[z_0,\cdots,z_{n+1}]\in {\PP}^{n+1}: -
\sum_{j=0}^{\ell}|z_j|^2 + \sum_{j=\ell+1}^{n+1}|z_j|^2 <0
\right\}.$$ The boundary $\p{{\BB}_\ell^{n+1}}$  is locally
holomorphically equivalent to the hyperquadric $ {\mathbb
H}^{n+1}_\ell\subset {\CC}^{n+1}$ of signature $\ell$ defined by
$\Im{z_{n+1}}=-\sum_{j=1}^{\ell}|z_j|^2+
\sum_{j=\ell+1}^{n+1}|z_j|^2,$ where  $(z_1,\cdots, z_{n+1})$ is the
coordinates of ${\CC}^{n+1}$.

For $0<\epsilon<<1$, ${M_\epsilon}$ is a compact smooth
real-algebraic hypersurface with Levi form non-degenerate of  the
same signature $\ell$.

\begin{theorem}\label{thm222} ([HZa])
There is an $\epsilon_0>0$ such that for
$0<\epsilon<\epsilon_0$, the following holds: (i)
 ${M_\epsilon}$ is a
smooth real-algebraic hypersurface in ${\PP}^{n+1}$ with
non-degenerate Levi form of signature $\ell$ at every point. (ii)
There does not exist  any holomorphic embedding from any open piece
of ${M_\ell}$ into $ {\mathbb H}_\ell^{N+1}$. \end{theorem}

When $0<\epsilon<<1$, since ${M_\epsilon}$ is a small algebraic
deformation of the generalzied sphere, we see that ${M_\epsilon}$ must also be a compact
real-algebraic Levi non-degenerate hypersurface in ${\PP}^{n+1}$
with signature $\ell$ diffeomorphic to the generalized sphere which is the boundary of the generalized  ball
${\BB}^{n+1}_\ell\subset {\PP}^{n+1}$.

\medskip
{\it Proof of Theorem \ref{thm222}}:
The proof uses the following algebraicity of the first author:

\begin{theorem}[Hu2, Corollary  in $\S 2.3.5$] \label{thm333}  Let
$M_1\subset {\CC }^n$ and $M_2\subset {\CC}^N$ with $N\ge n\ge 2$ be
two Levi non-degenerate real-algebraic hypersurfaces.
Let $p\in M_1$ and $U_p$ be a small connected open neighborhood of
$p$ in ${\mathbf C}^n$ and $F$ be a holomorphic map from $U_p$ into
${\mathbf C}^N$ such that $F(U_p\cap M_1)\subset M_2$ and
$F(U_p)\not \subset M_2$. Suppose that $M_1$ and $M_2$ have the same
signature $\ell$ at $p$ and $F(p)$, respectively. Then $F$ is
algebraic in the sense that each component of $F$ satisfies a
nontrivial holomorphic polynomial equation. \end{theorem}
Next,  we
compute the Chern-Moser-Weyl tensor of $M_\epsilon$ at the point
$$P_0:=[\xi^0_0,\cdots,\xi^0_{n+1}], \quad \xi^0_j=0 \text{ for } j\not = 0,\ell+1,
 \quad \xi^0_0=1,\quad \xi_{\ell +1}^0=1,$$
 and consider the coordinates
$$\xi_0=1,\quad \xi_j=\frac{\eta_j}{1+\sigma}, \quad j=1,\cdots, \ell,
 \quad \xi_{\ell+1}=\frac{1-\sigma}{1+\sigma},
\quad \xi_{j+1}=\frac{\eta_j}{1+\sigma},\quad j=\ell+1 ,\cdots, n.$$
Then in the $(\eta,\sigma)$-coordinates, $P_0$ becomes the origin
and $M_\epsilon$ is defined near the origin by an equation in the
form:

\begin{equation}\label{26}
\rho=-4\Re{\sigma}-\sum_{j=1}^{\ell}|\eta_j|^2+\sum_{j=\ell+1}^{n}|\eta_j|^2
+{a}(|\eta_1|^4-|\eta_{n}|^4)+o(|\eta|^4)=0,
\end{equation}
for some $a>0$. Now, let $Q(\eta,\-\eta)=-a(|\eta_1|^4-|\eta_n|^4)$
and make a standard $\ell$-harmonic decomposition [SW]:

\begin{equation}\label {33}
Q(\eta,\-\eta)=N^{(2,2)}(\eta,\-\eta)+
A^{(1,1)}(\eta,\-\eta)|\eta|^2_{\ell}.
\end{equation}
Here $N^{(2,2)}(\eta,\eta)$ is a $(2,2)$-homogeneous polynomial in
$(\eta,\-{\eta})$ such that $\Delta_\ell N^{(2,2)}(\eta,\-\eta)=0$
with $\Delta_\ell$ as before. Now $N^{(2,2)}$ is the
Chern-Moser-Weyl tensor of $M_\epsilon$ at $0$ (with respect to an
obvious contact form) with $N^{(2,2)}(\eta,\-\eta)=Q(\eta,\-\eta)$
for any $\eta\in{\mathcal C}T^{(1,0)}_0 M_e$. Now the value of the
Chern-Moser-Weyl tension has negative  and positive value at
$X_1=\frac{\p }{\p \eta_1}+\frac{\p }{\p \eta_{\ell+1}}|_0$ and
$X_2=\frac{\p }{\p \eta_{2}}+\frac{\p }{\p \eta_{n}}|_0$,
respectively. If $\ell>1$, then  both $X_1$ and $X_2$ are in
${\mathcal C}T^{(1,0)}_0 M_e$. We see that the Chern-Moser-Weyl
tensor can not be pseudo semi-definite near the origin in such a
coordinate system.

Next, suppose an open piece $U$ of ${M_\epsilon}$ can be
holomorphically and transversally embedded into the ${\mathbf
H}_\ell^{N+1}$ for $N>n$ by $F$. Then by the algebraicity result in
Theorem \ref{thm333}, $F$ is algebraic. Since the branching points
of $F$ and the points where $F$ is not defined (poles or points of
indeterminancy of $F$) are contained in a complex-algebraic variety
of codimension at most one, $F$ extends holomorphically along a
smooth curve  $\gamma$ starting from some point in $U$ and ending up
at some point $p^* (\approx 0)\in M_\ell$ in the
$(\eta,\sigma)$-space where the Chern-Moser-Weyl tensor of $M_\epsilon$
is not pseudo-semi-definite. By the uniqueness of real-analytic
functions, the extension of $F$ must also map an open piece of $p^*$
into  $ {\bf H}^{N+1}_\ell$. The extension is not totally
degenerate. By Theorem \ref{thm111}, we get a contradiction.
$\endpf$

\bigskip

\section{Non-embeddability of compact strongly psuedo-convex real algebraic hypersurfaces into  spheres}

As discussed in the previous sections,  spheres serve as the model
of strongly pseudoconvex real hypersurfaces where the
Chern-Moser-Weyl tensor vanishes. An immediate application of the
invariant property for the Chern-Moser-Weyl tensor is that very rare
strongly pseudoconvex real hypersurfaces can be biholomorphically
mapped to a unit sphere. Motivated by various embedding theorems in
geometries (Nash embedding, Remmert embedding theorems, etc), a
natural question to pursue in Several Complex Variables is to
determine when a real hypersurface in $\mathbb{C}^n$ can be
holomorphically embedded  into the unit sphere
$\mathbb{S}^{2N-1}=\{Z \in \mathbb{C}^N: ||Z||^2=1\}.$

By a  holomorphic  embedding of $M\subset {\mathbb C}^n$ into $M'
\subset {\mathbb C}^N$, we mean a holomorphic embedding of an open
neighborhood $U$ of $M$ into a neighborhood $U'$ of $M'$, sending
$M$ into $M'$.
We also say $M$ is locally holomorphically embeddable into $M'$ at $p
\in M$, if there is a neighborhood $V$ of $p$ and a holomorphic
embedding $F: V \rightarrow \mathbb{C}^{N}$ sending $M \cap V$ into
$M'$.

A real hypersurface holomorphically embeddable into a sphere is
necessarily strongly pseudoconvex and real-analytic. However, due to
results by Forstneri\'c [For1] (See a recent work [For2] for further result) and Faran [Fa], not every strongly
pseudoconvex real-analytic hypersurface can be embedded into a
sphere. Explicit examples of non-embeddable strongly pseudoconvex
real-analytic hypersurfaces constructed much later in [Za1]. Despite
a vast of literature devoted to the embeddability problem, the
following question remains an open question of long standing. Here
recall a smooth real hypersurface in an open subset $U$ of
$\mathbb{C}^n$ is called real-algebraic, if it has a real-valued
polynomial defining function.

\begin{question} \label{question} Is every compact real-algebraic strongly pseudoconvex real
hypersuraface in $\mathbb{C}^n$ holomorphically embeddable into a
sphere of sufficiently large dimension?
\end{question}

 Part of the motivation to study this embeddability problem is a
well-known result due to Webster [We2] which states that every
real-algebraic Levi-nondegenerate hypersurface admits a transversal
holomorphic embedding into a non-degenerate hyperquadric in
sufficiently large complex space. (See also [KX] for further study
along this line.) Notice  that in [HZa], the authors showed that
there are many compact real-algebraic pseudoconvex real
hypersurfaces with just one weakly pseudoconvex point satisfying the
following property:  Any open piece of them cannot be
holomorphically embedded into any compact real-algebraic strongly
pseudoconvex hypersurfaces which, in particular,  includes spheres.
Many other related results can be found in  the work of Ebenfelt-Son
[ES], Fornaess [Forn], etc.

In [HLX], the authors constructed the following family of compact
real-algebraic strongly pseudoconvex real hypersurfaces:
\begin{equation} \label{equhym}
M_{\epsilon}=\{(z,w)\in\mathbb{C}^2:
\varepsilon_0(|z|^8+c\mathrm{Re}|z|^2z^6)+|w|^2+ |z|^{10}+
{\epsilon} |z|^2-1=0\},~~ 0 < \epsilon < 1.
\end{equation}
Here, $2<c<\frac{16}{7}$, $\varepsilon_0 >0$ is a sufficiently small
number  such that $M_\varepsilon$ is smooth for all $0 \leq \epsilon
<1$. An easy computation shows that for any $0 < \epsilon < 1, M_{\epsilon}$ is strongly pseudoconvex.
$M_{\epsilon}$ is indeed a small algebraic deformation of the boundary of the
famous Kohn-Nirenberg domain [KN]. It is shown in [HLX] that for any
integer $N,$ there exists a small number $0 < \epsilon(N) < 1$, such
that for any $0 < \epsilon < \epsilon(N)$, $M_{\epsilon}$ cannot be
locally holomorphically embedded into the unit sphere
$\mathbb{S}^{2N-1}$ in $\mathbb{C}^N$. More precisely, any
holomorphic map sending an open piece of $M_{\epsilon}$ to
$\mathbb{S}^{2N-1}$ must be a constant map. We will write
$$\rho_{\epsilon}=\rho_{\epsilon}(z, w, \overline{z}, \overline{w}):=\varepsilon_0(|z|^8+c\mathrm{Re}|z|^2z^6)+|w|^2+
|z|^{10}+ {\epsilon} |z|^2-1.$$

  We first fix some notations. Let $M \subset
\mathbb{C}^{n}$ be a real-algebraic subset defined by a family of
real-valued polynomials $\{\rho_{\alpha}(Z,\overline{Z})=0\},$ where
$Z$ is the coordinates of $\mathbb{C}^{n}.$  Then the
complexification $\mathcal{M}$ of $M$ is the complex-algbraic subset
in $\mathbb{C}^n \times \mathbb{C}^n$ defined by
$\rho_{\alpha}(Z,W)=0$ for each $\alpha, (Z, W) \in \mathbb{C}^n
\times \mathbb{C}^n.$ Then for $p \in \mathbb{C}^{n},$ the Segre
variety of $M$ associated with the point $p$ is defined by
$Q_{p}:=\{Z \in \mathbb{C}^n:(Z,\overline{p}) \in \mathcal{M}\}.$
The geometry of Segre varieties of a real-analytic hypersurface has
been used in many literatures since the work of Segre [S] and
Webster [We].

In this note, fundamentally based on our previous joint work with Li
[HLX], we show that $M_{\epsilon}$ cannot be locally holomorphically
embedded into any unit sphere. The other important observation we
need is the fact that for some $p\in M_\epsilon$, the associated
Segre variety $Q_p$ cuts $M_\epsilon$ along a one dimensional real
analytic subvariety inside $M_\epsilon$. The geometry related to
intersection of  the Segre variety with the boundary plays an
important role in the study of many problems in Several Complex
Variables. (We mention, in particular, the work of  D'Angelo-Putinar
[DP], Huang-Zaitsev [HZa]).


 This then provides a counter-example to a long
standing open question--- Question \ref{question}. (See [HZa] for more
discussions on this matter).

\begin{theorem}\label{t0}
There exist compact real-algebraic strongly pseudoconvex
real hypersurfaces in $\mathbb{C}^2$, diffeomorphic to the sphere, that are not locally
holomorphically embeddable into any sphere. In particular,  for
sufficiently small positive $\varepsilon_0, \epsilon, M_{\epsilon}$ cannot be
locally holomorphically embedded into any sphere. More precisely, a local
holomorphic map sending an open piece of $M_{\epsilon}$ to a unit
sphere must be a  constant map.
\end{theorem}

 Write $D_{\epsilon}= \{ \rho_{\epsilon} < 0 \}$ as the
interior domain enclosed by $M_{\epsilon}.$ Since $M_{\epsilon}$ is
a small smooth deformation of $\{ |z|^{10} + |w|^2=1 \}$ for small
$\varepsilon_0$ and $\epsilon$. This imples $M_{\epsilon}$ is diffeomorphic to the unit sphere $\mathbb{S}^3$ for sufficiently small $\varepsilon_0$ and $\epsilon$.
Consequently, $M_{\epsilon}$  separates ${\CC}^2$
into two connected components $D_\epsilon$ and ${\CC }^2\sm
\overline{D_\epsilon}$.







\begin{proposition}\label{prop1e}
Let $p_0=(0,1) \in M_{\epsilon}.$ Let $Q_{p_0}$ be the Segre variety
of $M_{\epsilon}$ associated to $p_0.$ There exists
$\widetilde{\epsilon}>0$ such that for each $0 < \epsilon <
\widetilde{\epsilon}$,  $Q_{p_0} \cap M_{\epsilon}$ is a real
analytic subvariety of   dimension one.
\end{proposition}

{\it Proof of Proposition \ref{prop1e}}: It suffices to show that
there exists $q \in Q_{p_{0}}$ such that $q \in D_{\epsilon}.$ Note
that $ Q_{p_0}=\{(z, w): w=1 \}$. Set
$$\psi(z, \epsilon)=\varepsilon_0(|z|^8 + c\mathrm{Re} |z|^2 z^6)+|z|^{10}+ \epsilon |z|^2,~0\leq  \epsilon < 1.$$
Note $q=(\mu_0, 1) \in D_{\epsilon}$ if and only if $\psi(\mu_0, \epsilon)<
0.$   Now, set $\phi(\lambda, \epsilon)=\varepsilon_0 \lambda^8 (1
-c)+ \lambda^{10}+ \epsilon \lambda^2, 0 \leq  \epsilon < 1.$
 First we note there exists small $\lambda'> 0,$ such that $\phi(\lambda', 0)< 0$. Consequently, we can find $\widetilde{\epsilon} >0$ such that for each $0 < \epsilon \leq \widetilde{\epsilon}, \phi(\lambda', \epsilon) <0.$
Write $\mu_0=\lambda' e^{i\frac{\pi}{6}}.$ It is easily to see that
$\psi(\mu_0, \epsilon)< 0$ if $0 < \epsilon \leq \widetilde{\epsilon}$. This establishes Proposition \ref{prop1e}.
$\endpf$

\bigskip

\bigskip

\begin{proposition}\label{prpt01}
Let $M:= \{Z \in \mathbb{C}^n: \rho(Z, \overline{Z})=0\},n \geq 2, $
be a compact, connected, strongly pseudo-convex real-algebraic
hypersurface.  Assume that there exists a point $p \in M$ such that
the associated Segre variety $Q_{p}$ of $M$ is irreducible and
$Q_{p}$ intersects $M$ at infinitely many points. Let $F$ be a
holomorphic rational map sending an open piece of  $M$ to the unit
sphere $\mathbb{S}^{2N-1}$ in some $\mathbb{C}^N.$ Then $F$ is a
constant map.
\end{proposition}

{\it Proof of Proposition \ref{prpt01}}: Let $D$ be the interior domain enclosed by $M.$ From the assumption and a
theorem of Chiappari [Ch], we know $F$ is holomorphic  in a
neighborhood $U$ of $\overline{D}$ and sends
$M$ to $\mathbb{S}^{2N-1}.$ Consequently, if we write
$\mathcal{S}$ as the singular set of $F$, then it does not intersect  $U$. Write $Q'_q$ for the Segre variety of
$\mathbb{S}^{2N-1}$ associated to $q \in \mathbb{C}^N$.  We first conclude by complexification
that for a small neighborhood $V$ of $p,$
\begin{equation}\label{eqnpq}
F(Q_p \cap V) \subset Q'_{F(p)}.
\end{equation}
Note that $\mathcal{S} \cap Q_p$ is a Zariski close proper
subset  of $Q_p$. Notice that  $Q_p$ is connected as it is
irreducible. We conclude by unique continuation that if
$\widetilde{p} \in Q_p$ and $F$ is holomorphic at $\widetilde{p}$,
then $F(\widetilde{p}) \in Q'_{F(p)}$. In particular, if
$\widetilde{p} \in Q_p \cap M,$ then $F(\widetilde{p}) \in Q'_{F(p)}
\cap \mathbb{S}^{2N-1}=\{ F(p) \}.$ That is,
$F(\widetilde{p})=F(p).$

\smallskip

Notice by assumption that $Q_{p} \cap M$ is a compact set and contains infinitely
many points.  Let $\hat{p}$ be an accumulation point of $Q_{p} \cap
M.$ Clearly, by what we argued above, $F$ is not one-to-one in any
neighborhood of $\hat{p}.$  This shows that $F$ is constant. Indeed,
suppose $F$ is not a constant map. We then conclude that $F$ is a
holomorphic embedding near $\hat{p}$ by a standard Hopf lemma type argument (see [Hu2], for instance) for both
$M_\epsilon$ and ${\mathbb S}^{2N-1}$ are strongly pseudo-convex.
This completes the proof of Proposition \ref{prpt01}. $\endpf$

\medskip
{\it  Proof of Theorem \ref{t0}:} Pick $p_{0}=(0,1) \in
M_{\epsilon}.$   Notice that the associated Segre variety $Q_{p_0}=\{(z,1): z
\in \mathbb{C} \}$ is an irreducible complex variety in
$\mathbb{C}^2$.  Let $\epsilon, \varepsilon_0$ be
sufficiently small such that Proposition \ref{prop1e}  holds.

Now,  let $F$ be a holomorphic  map defined in a
small neighborhood $U$ of  some point $q \in M_{\epsilon}$ that sends an open
piece of $M_{\epsilon}$ into  $\mathbb{S}^{2N-1},  N \in \mathbb{N}$. It is shown in
[HLX] that $F$  is a rational map. Then it follows from Proposition \ref{prpt01} that
$F$ is a constant map. We have thus established Theorem \ref{t0}.
$\endpf$

\end{document}